\documentclass{gtart_a}
\pdfoutput=1

\usepackage{xy,graphicx}
\xyoption{all}

\title[Moment mappings on loop groups]{Connectivity properties of 
moment maps\\on based loop groups}

\author[Harada]{Megumi Harada}
\givenname{Megumi}
\surname{Harada}
\address{Department of Mathematics and Statistics\\
McMaster University\\1280 Main Street West\\\newline 
Hamilton\\Ontario L8S 4K1\\Canada}
\email{Megumi.Harada@math.mcmaster.ca}
\urladdr{}

\author[Holm]{Tara S Holm}
\givenname{Tara S}
\surname{Holm}
\address{Department of Mathematics\\589 Malott Hall\\
Cornell University\\\newline
Ithaca\\NY 14850-4201\\USA}
\email{tsh@math.cornell.edu}
\urladdr{}

\author[Jeffrey]{Lisa C Jeffrey}
\givenname{Lisa C}
\surname{Jeffrey}
\address{Department of Mathematics\\University of Toronto\\
Toronto\\Ontario M5S 2E4\\Canada}
\email{jeffrey@math.toronto.edu}
\urladdr{}

\author[Mare]{Augustin-Liviu Mare}
\givenname{Augustin-Liviu}
\surname{Mare}
\address{Department of Mathematics and Statistics\\
University of Regina\\College West 307.14\\\newline
Regina\\Saskatchewan S4S 0A2\\Canada}
\email{mareal@math.uregina.ca}
\urladdr{}

\volumenumber{10}
\issuenumber{}
\publicationyear{2006}
\papernumber{40}
\lognumber{0589}
\startpage{1607}
\endpage{1634}

\doi{}
\MR{}
\Zbl{}

\keyword{}
\subject{primary}{msc2000}{}
\subject{secondary}{msc2000}{}

\received{4 April 2005}
\revised{15 July 2005}
\accepted{6 September 2006}
\proposed{Ralph Cohen}
\seconded{Haynes Miller, Frances Kirwan}
\published{28 October 2006}
\publishedonline{28 October 2006}
\corresponding{}
\editor{CPR}
\version{}

\arxivreference{math.SG/0503684}

\let\xysavmatrix\xymatrix
\def\xymatrix{\disablesubscriptcorrection\xysavmatrix}
\AtBeginDocument{\let\bar\wbar\let\tilde\wtilde\def\bC{\mathbb{C}}
\def\C{C_{\lambda}}\def\t{\mathfrak{t}}\def\P{\mathbb{P}}
\def\S{Section }\def\notin{\not\in}}

\makeatletter
\def\cnewtheorem#1[#2]#3{\newtheorem{#1}{#3}[section]
\expandafter\let\csname c@#1\endcsname\c@theorem}

\numberwithin{equation}{section}

\newtheorem{theorem}{Theorem}[section]
\cnewtheorem{corollary}[theorem]{Corollary}
\cnewtheorem{proposition}[theorem]{Proposition}
\cnewtheorem{lemma}[theorem]{Lemma}

\theoremstyle{definition}
\cnewtheorem{definition}[theorem]{Definition}
\cnewtheorem{remark}[theorem]{Remark}

\makeatother  

\def\g{\mathfrak{g}}

\def\g{\frak{g}}
\def\a{\rm alg}
\def\t{\mathfrak{t}}

\newcommand{\into}{{\hookrightarrow}}

\begin{document}

\begin{asciiabstract}
For a compact, connected, simply-connected Lie group G, the loop group
LG is the infinite-dimensional Hilbert Lie group consisting of
H^1-Sobolev maps S^1-->G.  The geometry of LG and its homogeneous
spaces is related to representation theory and has been extensively
studied.  The space of based loops Omega(G) is an example of a
homogeneous space of $LG$ and has a natural Hamiltonian T x S^1
action, where T is the maximal torus of G.  We study the moment map mu
for this action, and in particular prove that its regular level sets
are connected. This result is as an infinite-dimensional analogue of a
theorem of Atiyah that states that the preimage of a moment map for a
Hamiltonian torus action on a compact symplectic manifold is
connected. In the finite-dimensional case, this connectivity result is
used to prove that the image of the moment map for a compact
Hamiltonian T-space is convex. Thus our theorem can also be viewed as
a companion result to a theorem of Atiyah and Pressley, which states
that the image mu(Omega(G)) is convex.  We also show that for the
energy functional E, which is the moment map for the S^1 rotation
action, each non-empty preimage is connected.
\end{asciiabstract}

\begin{htmlabstract}
For a compact, connected, simply-connected Lie group G, the loop group
LG is the infinite-dimensional Hilbert Lie group consisting of
H<sup>1</sup>&ndash;Sobolev maps S<sup>1</sup> &rarr; G. The geometry
of LG and its homogeneous spaces is related to representation theory
and has been extensively studied.  The space of based loops
&Omega;(G) is an example of a homogeneous space of LG and has a
natural Hamiltonian T &times; S<sup>1</sup> action, where T is the
maximal torus of G.  We study the moment map &micro; for this action,
and in particular prove that its regular level sets are
connected. This result is as an infinite-dimensional analogue of a
theorem of Atiyah that states that the preimage of a moment map for a
Hamiltonian torus action on a compact symplectic manifold is
connected. In the finite-dimensional case, this connectivity result is
used to prove that the image of the moment map for a compact
Hamiltonian T&ndash;space is convex. Thus our theorem can also be
viewed as a companion result to a theorem of Atiyah and Pressley,
which states that the image &micro;(&Omega;(G)) is convex.  We also
show that for the energy functional E, which is the moment map for the
S<sup>1</sup> rotation action, each non-empty preimage is connected.
\end{htmlabstract}

\begin{abstract}
For a compact, connected, simply-connected Lie group $G$, the loop
group $LG$ is the infinite-dimensional Hilbert Lie group consisting of
$H^1$--Sobolev maps \(S^1 \to G.\) The geometry of $LG$ and its
homogeneous spaces is related to representation theory and has been
extensively studied.  The space of based loops \(\Omega(G)\) is an
example of a homogeneous space of $LG$ and has a natural Hamiltonian
$T \times S^1$ action, where $T$ is the maximal torus of $G$.  We
study the moment map $\mu$ for this action, and in particular prove
that its regular level sets are connected. This result is as an
infinite-dimensional analogue of a theorem of Atiyah that states that
the preimage of a moment map for a Hamiltonian torus action on a
compact symplectic manifold is connected. In the finite-dimensional
case, this connectivity result is used to prove that the image of the
moment map for a compact Hamiltonian $T$--space is convex. Thus our
theorem can also be viewed as a companion result to a theorem of
Atiyah and Pressley, which states that the image $\mu(\Omega(G))$ is
convex.  We also show that for the energy functional $E$, which is the
moment map for the $S^1$ rotation action, each non-empty preimage is
connected.
\end{abstract}

\maketitle

\section{Introduction}\label{secintro}

The main results of this paper are infinite-dimensional analogues of
well-known results in finite-dimensional symplectic geometry.  More
specifically, given a compact connected Hamiltonian $T$--manifold with
moment map \(\mu\co  M \to \t^*,\) a theorem of Atiyah \cite[Theorem 1]{At} states that
any non-empty preimage of $\mu$ is connected. 
This connectivity result is intimately
related to the famous convexity result of Atiyah and
Guillemin--Sternberg \cite{At,Gu-St}, that states that the image
\(\mu(M) \subseteq \t^*\) is {\em convex.} Indeed, in the original
paper of Atiyah, the connectivity result is used to establish the
convexity of $\mu(M)$.

\begin{figure}[ht!]
\begin{center}
\includegraphics[width=1.6in]{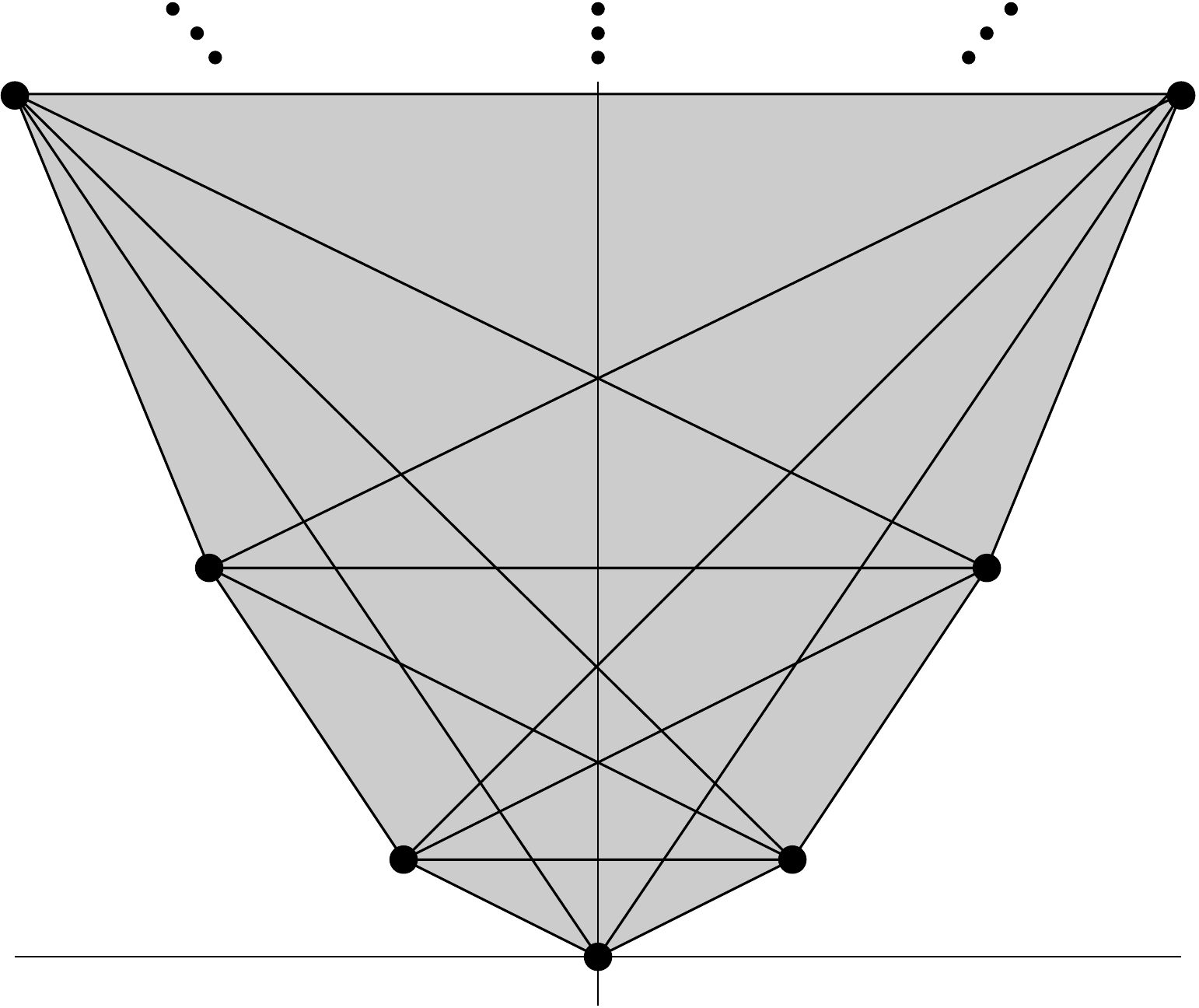}
\end{center}
\caption{ The shaded region indicates a portion of the image in
$(\t^2)^*$ of the moment map for the $T^2 = T^1 \times S^1$ action on
$\Omega(SU(2))$, which is the convex hull of the integer points on a
parabola.  (We have also drawn in the critical values of $\mu$, which
are the line segments.)  The convexity of the image was proven in
[26].  We will show that the level sets of this map are connected.}
\label{figOmegaSUtw}
\end{figure}

Atiyah and Pressley \cite{At-Pr} showed that the convexity results
mentioned above generalize to a certain class of infinite-dimensional
manifolds that are homogeneous spaces of loop groups $LG$.  An example
is illustrated in \fullref{figOmegaSUtw}. Indeed, there is a general
principle (see eg, Pressley and Segal \cite{Pr,Pr-Se}) which says that
these infinite-dimensional homogeneous spaces of $LG$ behave, in many
respects, like {\em compact} K\"ahler manifolds. The result of
\cite{At-Pr} may be viewed as an instance of this, as should the
results of this paper.  However, in contrast to the argument given in
\cite{At}, the convexity result in \cite{At-Pr} is proven without any
mention of the connectivity of level sets of the moment map. Therefore
a natural question still remains: is there an analogous connectivity
result for these homogeneous spaces of loop groups? The purpose of
this paper is to answer this question in the affirmative.

In the finite-dimensional setting, connectivity of the level sets has
several important implications.  Some of these suggest very natural
questions for the infinite-dimensional context; we mention three of
these here. First, connectivity along with a local normal form
description of the level sets yields the Atiyah--Guillemin--Sternberg
convexity result.  A new proof of the Atiyah--Pressley convexity result
for $\Omega G$ that uses our connectivity result and some additional
information about the local normal forms would be very
interesting. Second, connectivity is the first piece of Kirwan's
surjectivity theorem \cite{Ki1}.  Given a compact connected
Hamiltonian $G$--manifold with moment map \(\mu\co M \to \g^*,\) Kirwan
showed that the inclusion of a level set $\mu^{-1}(0)\hookrightarrow
M$ induces a surjection in equivariant cohomology.  If the level set
is regular, its equivariant cohomology is isomorphic to the ordinary
cohomology of the symplectic quotient.  Thus, there is a surjection
$H_G^*(M) \to H^*(M/\! /G)$.  For the homogeneous spaces of loop
groups, connectivity of the level sets implies that this map is a
surjection for degree $0$ cohomology.  A surjection $\kappa\co H^*_{T
\times S^1}(\Omega G) \to H^*(\Omega G /\!/ T\times S^1)$ would be a
very powerful result. By understanding the map $\kappa$ and its
kernel, and using the description of $H^*_{T \times S^1}(\Omega G)$ in
Harada, Henriques and Holm \cite{Ha-He-Ho}, one could obtain a
description of the cohomology ring of this quotient.  Extending known
results about $\kappa$ and its kernel (eg, Jeffrey, Kirwan, Tolman and
Weitsman \cite{JK:LocalizationNonabelian,TW:symplecticquotients}) will
nevertheless require significant technical prowess.  Moreover, there
is a possible interpretation in terms of representation
theory. According to the so-called ``quantization commutes with
reduction'' principle in the finite-dimensional setting, it may be
possible to interpret the spaces $\Omega G /\!/ T \times S^1$ as
geometric analogues of $T\times S^1$--weight spaces for certain loop
group representations.  Third, in finite dimensions, Karshon and
Lerman \cite{Ka-Le} used the connectivity of level sets to deal with
the double commutator conjecture of Guillemin and Sternberg
\cite{Gu-St-84}. This conjecture states that the centralizer of the
algebra of $G$--invariant functions on a Hamiltonian $G$--manifold is
{\em collective}, ie,\ the pullback via the $G$--moment map of a
smooth function on $\g^*$.  The connectivity results contained herein
may aid in the determining the validity of an infinite dimensional
analogue of this theorem.

We now state more precisely the main results of this paper.  Let $G$ be a
compact, connected and simply-connected Lie group. We fix a
$G$--invariant inner product $\langle \ , \ \rangle$ on the Lie algebra
$\g$.  Let \(L_1(G)\) denote the space of all maps \(S^1 \cong
\R/2\pi\Z \to G\) that are of Sobolev class $H^1$, ie,
\[
L_1(G) := \{ \eta \in H^1(S^1, G) \}.
\]

{\bf Remark}\qua Throughout the remainder
of this paper, $H^1$ stands always for 
``Sobolev class" (and not for ``cohomology"). 

The space $L_1(G)$ is a group by pointwise multiplication. We now
consider the subgroup $\Omega_1(G)$ of $L_1(G)$ consisting of {\em
based} loops $\gamma$, ie, we add the requirement \(\gamma(0) = e \in
G.\) The group $L_1(G)$ acts by conjugation on $\Omega_1(G)$ and it is
straightforward to see that the stabilizer of the constant loop at $e$
is exactly the set \(G \subseteq L_1(G)\) of constant loops. Hence
$\Omega_1(G)$ is a homogeneous space \(L_1(G)/G.\) This is a smooth
Hilbert manifold (see eg, Palais \cite{Pa}) which is also K\"ahler
(see Atiyah, Pressley and Segal \cite{At-Pr,Pr-Se}), and is the space
on which we will focus for the remainder of the paper.

We consider the following torus action on $\Omega_1(G)$. The
maximal torus $T$ of $G$ is a subgroup of $L_1(G)$, and hence acts on
the left on $\Omega_1(G) \cong L_1(G)/G$. More specifically, for \(t
\in T, \gamma \in \Omega_1(G),\) the action is defined by
\[
(t \gamma)(\theta) := t \gamma(\theta) t^{-1}.
\]
In addition to this maximal torus action, we also have a {\em rotation
action} of $S^1$ that rotates the loop variable. For \(e^{is} \in
S^1\) and \(\gamma \in \Omega_1(G),\) we have 
\begin{equation}\label{rotation}
(e^{is} \gamma)(\theta) := \gamma(s + \theta)\gamma(s)^{-1}.
\end{equation}
These actions commute, so we obtain a \(T \times S^1\) action on
$\Omega_1(G)$ which is also Hamiltonian with respect to the K\"ahler
structure on $\Omega_1(G)$. The $S^1$--moment map is in fact a
well-known functional on spaces of loops: it is the {\em energy
functional}
\begin{equation}\label{energyon}
E(\gamma) := \frac{1}{4\pi} \int_0^{2\pi} \|\gamma(\theta)^{-1}
\gamma'(\theta) \|^2 d\theta.
\end{equation}
The moment map for the $T$--action is given by a similar functional,
\begin{equation}\label{eqTmoment}
p(\gamma) = {\rm pr}_{\t}\left(
\frac{1}{2\pi} \int_0^{2\pi} \gamma(\theta)^{-1}\gamma'(\theta) d\theta\right),
\end{equation}
where $\t:={\rm Lie}(T)$ and ${\rm pr}_{\t} \co \g \to \t$ denotes the
orthogonal projection with respect to the fixed $G$--invariant inner
product $\langle \cdot , \cdot \rangle$. Using this inner product, we
may identify \(\t^*\) with $\t$, and use the standard inner product on
$\R$ to identify $\R^*$ with $\R$. With these identifications, the
$T\times S^1$--moment map \(\mu\co  \Omega_1(G) \to \t^* \oplus \R^* \cong
\t \oplus \R\) is given by \(\mu = p \oplus E.\)

There is a subspace of $\Omega_1(G)$ that
we will also study, namely  $\Omega_{\a}(G)$, which is the space of algebraic  loops
in $G$. The main result of \cite{At-Pr}, as mentioned above,
is that the images \(\mu(\Omega_{\a}(G))\) and \(\mu(\Omega_1(G))\) 
are both convex. 
Our results are as follows. 

\begin{theorem}\label{firstmain} Any level of $\mu\co  \Omega_{\a}(G) \to \t \oplus \R$ is connected or empty.
\end{theorem}

\begin{theorem}\label{secondmain} Any regular level of $\mu\co \Omega_1 (G) \to \t \oplus \R$ is connected or empty.
\end{theorem}

We also have a statement of connectivity for the level sets of just
the energy functional, considered as a function on $\Omega_1(G)$ or on 
$\Omega_{\a}(G)$. In this case, we have the
connectivity for any (possibly
singular) value of $E$.

\begin{theorem}\label{zeromain}
For $\bullet \in \{1, \a\}$, all preimages of $E\co  \Omega_{\bullet}(G)\to \R$ are
either empty or connected.
\end{theorem}

We now make explicit the topologies on the spaces $\Omega_{\bullet}(G)$
with respect to which we state our
Theorems~\ref{firstmain}--\ref{zeromain} above.  First, $\Omega_1(G)$
 has a natural topology induced from
$H^1(S^1,\g)$ via the exponential
map $\exp \co  \g \to G$, which is a local homeomorphism around $0$ (see
\cite[Definition 2.4]{At-Pr}).    The issue of the topology on the subspace $\Omega_{\a}(G)$
is more subtle.  On the one hand, $\Omega_{\a}(G)$ is
equipped with the subset topology induced from the inclusion
\(\Omega_{\a}(G) \into \Omega_1(G).\) However, $\Omega_{\a}(G)$ also
has a stratification by algebraic varieties, which gives it a direct
limit topology (see \cite[\S 3.5]{Pr-Se}).  For Theorems~\ref{firstmain}
and~\ref{zeromain} above, we consider $\Omega_{\a}(G)$ with the {\em
direct limit} topology.  The two 
topologies on $\Omega_{\a}(G)$ are related; the subset topology is
 coarser than the direct limit topology (see eg, \fullref{finer}),
so \fullref{firstmain} in fact implies that $\mu^{-1}(a)\cap
\Omega_{\a}(G)$ is connected with respect to either of the possible
topologies of $\Omega_{\a}(G)$.

It is worth remarking here that the proofs of connectivity for
$\Omega_{\a}(G)$ and $\Omega_1(G)$ are markedly different in
flavor. To prove connectivity of levels for both $\mu$ and $E$ on
$\Omega_{\a}(G)$, we exploit the algebraic structure of
$\Omega_{\a}(G)$, in particular that there it may be described as a
union of finite-dimensional subvarieties. The proofs of the connectivity results for
$\Omega_1(G)$, on the other hand, use the connectivity results for
$\Omega_{\a}(G)$ and the density (in the $H^1$ topology) of
$\Omega_{\a}(G)$ as a subset of $\Omega_1(G)$. The tool used here is
Morse theory in infinite dimensions. To use Morse theory in this
context requires additional
technical hypotheses, such as the Palais--Smale condition (C). These
hypotheses must explicitly be checked in order for us to use the
Morse-theoretic arguments, and these are the main technical
difficulties in this paper.

We now outline the contents of this manuscript. In
\fullref{secbackground} we briefly recall important known results
about $\Omega_{\bullet}(G)$
 and its moment maps. Then in
\fullref{secalgebraic} we prove the result for algebraic loops
$\Omega_{\a}(G)$. \fullref{secregular} is devoted to the proof of
\fullref{secondmain}. The argument in Sections~\ref{secalgebraic}
and~\ref{secregular} also proves certain cases of
\fullref{zeromain}. Finally, in \fullref{secsingular}, we
prove \fullref{zeromain} for the case of the singular levels of
the energy functional on $\Omega_1(G)$, which requires a
separate argument.

\medskip
{\bf Remark}\qua It would be interesting to extend Theorem 1.1
and Theorem 1.2 to arbitrary adjoint orbits of Kac-Moody groups or
even to isoparametric submanifolds in Hilbert space (cf Terng and
Mare \cite{Te2,Ma}).

\medskip
{\bf Acknowledgments}\qua We would like to thank M Brion, R Cohen, Y
Karshon, Y-H Kiem, E Lerman, R Sjamaar and J Woolf for helpful
discussions.  We also thank the American Institute of Mathematics for
hosting all four authors while some of this work was conducted. The
second author was supported in part by a National Science Foundation
Postdoctoral Fellowship. The third and fourth authors are supported in
part by NSERC.

\section{Background material}\label{secbackground}

In this section we collect facts that will be needed in the proofs of
our main results.

\subsection{Loop groups}\label{subsecalgebraicloops}

The main reference for this section is \cite{At-Pr}, especially
section 2 (see also Freed \cite[\S 1]{Fr} or Pressley \cite[\S
2]{Pr}).  By definition, $L_1(G)=H^1(S^1, G)$ is the space of free
loops in $G$ of Sobolev class $H^1$.  These are maps $\eta\co S^1 \to
G$ with the property that for any local coordinate system $\Phi$ on
$G$, the map $\Phi \circ \eta \co S^1 \to \R^{\dim G}$ is of Sobolev
class $H^1$.  The space $L_1(G)$ is an infinite dimensional Lie group
(cf \cite[\S 13]{Pa}).  It carries a natural left invariant
Riemannian metric, which will be addressed as the $H^1$ metric. This
is uniquely determined by its restriction to the Lie algebra of
$L_1(G)$, which is $H^1(S^1, \g)$ (the tangent space at the constant
loop $e$).  First we fix an ${\rm Ad}(G)$--invariant metric on the Lie
algebra $\g$, denoted by $( \ , \ )$.  The $H^1$ metric is determined
by
$$\langle \gamma,\eta \rangle_e =\frac{1}{2\pi} \int_0^{2\pi} (\gamma(\theta), \eta(\theta))d\theta + 
 \frac{1}{2\pi} \int_0^{2\pi} (\gamma'(\theta), \eta'(\theta))d\theta,$$
 for $\gamma,\eta \in H^1(S^1,\g)$. 

As mentioned in the introduction, the object of study of our paper is
the space $\Omega_1(G)$. This is a closed submanifold of $L_1(G)$ and the $H^1$ metric defined above
induces a metric on $\Omega_1(G)$, which we will denote by the same
symbol $\langle \ , \ \rangle$ (the standard reference for this is
 \cite[\S 13, especially Theorem (6)]{Pa}).
 We will also consider the
K\"ahler metric on $\Omega_1(G)$. The details of the construction of
this metric can be found in the references indicated at the beginning
of the section (see also \cite[\S 8.9]{Pr-Se}).  We will just mention
here that $\Omega_1(G)$ carries a natural symplectic form $\omega$,
which is $L_1(G)$--invariant and its value at $e$ (the constant loop at the 
identity) is
$$\omega_e(\gamma,\eta) = \frac{1}{2\pi}\int_0^{2\pi} (\gamma'(\theta),\eta(\theta))d\theta.$$
There is a certain complex structure $J$ on $\Omega(G)$ and it can be shown that the triple
$(\Omega_1(G), \omega, J)$ is a K\"ahler (Hilbert) manifold
(see \cite[Proposition 8.9.8]{Pr-Se}).
 The corresponding 
K\"ahler metric 
will be denoted by $g$ in our paper. The difference between the $H^1$ and the K\"ahler metric
is that the first one is complete  whereas the second one is not (this observation will play an important role  in \fullref{secregular}). 
To understand where the difference  comes from, we note that both metrics are induced by
the closed embedding $\Omega_1(G)\subset L_1(G)$ (see above).
The $H^1$ metric on $L_1(G)=H^1(S^1,G)$ is obviously complete.
The K\"ahler metric on $\Omega_1(G)$ is the restriction of the $H^{\frac{1}{2}}$ metric
on $L_1(G)$ (this was first noted by Pressley in \cite[\S 2]{Pr}), hence it is not complete. For a detailed discussion about Sobolev metrics on loop groups, we refer the reader to \cite[\S 1]{Fr}.

Without loss of generality, we may assume that  $G$
  is a closed subgroup of $SU(N)$,  for $N$
sufficiently large. Moreover, we may also assume that the  torus
$T$  consists of diagonal matrices in $SU(N)$. It can be
seen \cite[\S 13]{Pa} that
$\Omega_1(G)$ is the space of all maps $\gamma$ of class $H^1$
from $S^1=\R/2\pi \Z$ to the space $M^{N\times N}(\bC)$ of all
complex $N\times N$ matrices that have  the properties that
$\gamma(S^1)\subset G$ and $\gamma(0)$ is the identity matrix $I$.
Morover, $\Omega_1(G)$ is a submanifold of $$L_1(M^{N\times N}(\bC)):=H^1(S^1, M^{N\times N}(\bC)).$$
Any $\gamma \in \Omega_1(G)$ has a Fourier expansion with
coefficients in  $M^{N\times N}(\bC)$. Then the subspace
$\Omega_{\a}(G)$ is defined to be the set of all algebraic loops,
ie, loops $\gamma$ with finite
Fourier series
\begin{equation}\label{ffinite} \gamma(\theta) = \sum_{k=-n}^n A_ke^{ik\theta},\end{equation} 
where $A_k\in M^{N\times N}(\bC)$, and $n\in \Z$. 
Consider the space
$$L_{\a}(M^{N\times N}(\bC)):=\{\gamma \co  S^1 \to M^{N\times N}(\bC)\ :
\ \gamma \ {\rm is \ a\ finite \ Fourier \ series} \}. $$
For any integer $n\ge 0$, we denote by $L_n$ the set of all elements 
of $L_{\a}(M^{N\times N}(\bC))$ of the form \eqref{ffinite}. The space $L_n$ can be
identified naturally with the direct product $(M^{N\times N}(\bC))^{2n+1}$, hence it carries a
natural topology. It is obvious that the spaces $L_0, L_1, \ldots $ are a filtration of 
$L_{\a}(M^{N\times N}(\bC))$, in the sense that
$$L_0 \subset L_1 \subset \ldots \subset L_n \subset \ldots L_{\a}(M^{N\times N}(\bC)),$$
$$ \bigcup_{j\ge 0} L_j = L_{\a}(M^{N\times N}(\bC)).\leqno{\hbox{and}} $$ 
In this way,   $L_{\a}(M^{N\times N}(\bC))$ can be equipped with the direct limit topology
(by definition, a subset $U\subset L_{\a}(M^{N\times N}(\bC))$ is open iff $U\cap L_j$ is open in $L_j$ for any $j\ge 0$).
The topology induced on 
 $\Omega_{\a}(G) = \Omega_1(G) \cap L_{\a} (M^{N\times N}(\bC))$ 
will be  also called the {\it direct limit topology}. Now $\Omega_{\a}(G)$ inherits a topology
from $\Omega_1(G)$ as well. This will be addressed as the {\it subspace topology}.
 The following result seems to be known
(see eg, \cite{Gu-Pr}).  We include a proof for the sake
of completeness (we will use this result in \fullref{secregular}) .

\begin{proposition}\label{finer} The direct limit topology on $\Omega_{\a}(G)$ is finer
than the subspace topology.
\end{proposition}
\begin{proof} In order to prove the lemma, it is sufficient to prove that the
direct limit topology is finer than (the restriction of the) $H^1$
topology on $L_{alg}(M^{N\times N}(\bC)).$
  Thus we need to prove that if
$\gamma_n \in L_{\a}(M^{N\times N}(\bC))$ and $(\gamma_n)$ converges
to 0 in the direct limit topology, then $(\gamma_n)$ converges to 0 in
the $H^1$ topology. Write
$$\gamma_n(\theta) =\sum_{k\in \Z} A_{kn}e^{ik\theta},$$ 
where $A_{kn} \in M^{N\times N}(\bC)$ and only finitely many \(A_{kn}
\neq 0\) for each $n$. 
Fix $K$ a positive integer. 
Consider $$V:=\{\gamma\in L_{\a}(M^{N\times N}(\bC)) {\rm  \ of \ type \ \eqref{ffinite}} \ : \
 |A_{k}| < 2^{-K-|k|} \  \forall k\in \Z\},$$ which is open in the direct limit topology and 
contains the matrix $0$.    The convergence of $(\gamma_n)$ to 0 in the direct limit topology
implies that there exists a positive integer $n(K)$ such that for any $n\ge n(K)$ we have
$\gamma_n\in V$, which means 
$$|A_{kn}| < 2^{-K-|k|}, \ \forall k\in \Z  .$$
In order to prove the convergence in the $H^1$ topology, recall that 
\[
\|\gamma_n\|_{H^1} := \frac{1}{2\pi}\int_0^{2\pi}(
\gamma_n(\theta), \gamma_n(\theta) ) + \frac{1}{2\pi}
\int_0^{2\pi} ( \gamma'_n(\theta), \gamma'_n(\theta)) =\sum_{k\in \Z} (1+ k^2)|A_{kn}|^2,
\]
which is less than
\[
 2^{-K} \sum_{k \in \Z} \frac{1+k^2}{2^{|k|}}.
\]
Because $\sum_{k \in \Z} \frac{1+k^2}{2^{|k|}}$ is finite, this 
implies $\|\gamma_n\|_{H^1}\to 0$, and hence $\gamma_n$ converges
to $0$ in the $H^1$ topology as desired. 
\end{proof}    

{\bf Remark}\qua In general, the direct limit topology on $\Omega_{\a}(G)$ is 
{\it strictly} finer than the subspace topology. For example, for $G=SU(2)$, we can consider
the sequence $\gamma_n\in \Omega_{\a}(SU(2))$ given by
\begin{eqnarray*}\gamma_n(z)  =  \left(%
\begin{array}{ccccccc}
  \sqrt{1-\frac{1}{n^4}} &  \frac{1}{n^2}z^n\\
  -\frac{1}{n^2}z^{-n} & \sqrt{1-\frac{1}{n^4}} \\
\end{array}%
\right)  
\end{eqnarray*}
where $z=e^{i\theta}$. 
We can see that $(\gamma_n)$ converges to the constant loop $I_2$ in the
$H^1$ topology, but not in the direct limit topology. To prove the first claim,
we note that
$$\left\|\gamma_n - I_2\right\|^2_{H_1} = \left\|\frac{1}{n^2}z^{-n}\right\|^2_{H_1} +
 2 \left(1-\sqrt{1-\frac{1}{n^4}}\right)^2
+\left\|\frac{1}{n^2}z^n\right\|^2_{H_1} $$
which is convergent to 0.  
In order to prove the nonconvergence of $(\gamma_n)$ in the direct limit topology,
we write
$$\gamma_n(z) = -\frac{1}{n^2}\left(%
\begin{array}{ccccccc}
  0 &  0\\
  1 & 0 \\
\end{array}%
\right) z^{-n} + \left(%
\begin{array}{ccccccc}
  \sqrt{1-\frac{1}{n^4}} &  0\\
  0 & \sqrt{1-\frac{1}{n^4}} \\
\end{array}%
\right) +  \frac{1}{n^2} \left(%
\begin{array}{ccccccc}
 0 &  1\\
 0 & 0\\
\end{array}%
\right)z^n .$$
Let $U$ be the (open) subspace of $L_{\a}(M^{2\times 2}(\bC))$
consisting of all series of the type \eqref{ffinite} where $n$ is arbitrary and  
$|A_k| <2^{-|k|}$, for all $k\in \Z$. 
Suppose that $(\gamma_n)$ converges to $I_2$ in the direct limit topology. This implies
that there exists $n_0$ such that for any $n\ge n_0$, we have $\gamma_n \in U$.
This implies
$$\left | \frac{1}{n^2} \left(%
\begin{array}{ccccccc}
  0 &  1\\
 0 & 0\\
\end{array}%
\right) \right| <2^{-n},$$
for all $n\ge n_0$, which is false.

\subsection{The Grassmannian model}\label{subsecloopgroups}
The main reference here is \cite{Pr-Se},  sections 7 and 8 (see also \cite[\S 2]{At-Pr}).
Let $H$ be the Hilbert space $L^2(S^1, \bC^N)$ and  $H_+$ the
(closed) subspace of $H$ consisting of all elements  with Fourier expansions of
the type
$$\sum_{k\ge 0}a_k e^{ik\theta},$$ where $a_k\in \bC^N$.
Denote by $H_-$ the orthogonal complement of $H_+$ in $H$. The
infinite-dimensional Grassmannian
$Gr(H)$ is the space of all closed linear subspaces $W \subset H$
such that 
\begin{itemize}
\item[(i)] the orthogonal projection ${\rm pr}_+\co W \to H_+$ is a Fredholm operator,
\item[(ii)]  the orthogonal projection ${\rm pr}_-\co W \to H_-$ is a Hilbert--Schmidt operator.
\end{itemize}
 This
Grassmannian $Gr(H)$ is a  Hilbert manifold which has a K\"ahler form
$\omega$ (see \cite[\S 7.8]{Pr-Se}). We also identify the following 
 important submanifolds of $Gr(H)$, that we will use in
the sequel. First, $Gr_{\infty}(H)$ (called the smooth Grassmannian) is the space  of all  $W\in Gr(H)$ for which  the images of both  orthogonal projections 
$W\to H_-$ and $W^{\perp}\to H_+$ consist of smooth functions.  By $Gr_{0}(H)$ one denotes  the space of all
 $W\in Gr(H)$ with the property that there exists an
integer $n\ge 0$ such that
\begin{equation}\label{H} e^{in\theta}H_+ \subset W \subset e^{-in\theta}H_+.\end{equation}
One can show that $Gr_0(H)\subset Gr_{\infty}(H)$ (see \cite[\S 7.2]{Pr-Se}). 

Like in section 2.1, we assume that $G$ is a subgroup of $SU(N)$.
The (linear) action of $T$ on $\bC^N$ induces in a natural way actions on $H$
and  on $Gr(H)$. 
 There also exists a
``rotation" action of $S^1$ on $H$, given by
$$(e^{i\theta}f)(z):=f(e^{i\theta}z),$$
for all $e^{i\theta},z\in S^1$ and all $f\in H$. It turns out that if $W\in Gr(H)$, then the space
$$e^{i\theta}W:=\{e^{i\theta}f  \ : \ f\in W\}$$
is also in $Gr(H)$. 
The actions of $T$ and $S^1$ on $Gr(H)$ commute with each other and induce an action of $T\times S^1$ on $Gr(H)$. The space $Gr_{\infty}(H)$ 
is $T\times S^1$ invariant (see \cite[\S 7.6]{Pr-Se}).  
Now if $W\in Gr(H)$, the map $S^1\to Gr(H)$ given by
$e^{i\theta}\mapsto e^{i\theta}W$ is in general not smooth.
The latter map is smooth provided that $W$ is in the smooth Grassmannian 
$Gr_{\infty}(H)$.   In this way, we 
obtain a smooth action of $T\times S^1$ on $Gr_{\infty}(H)$.  This
fits well with the symplectic structure on $Gr(H)$, as indicated in
the following proposition.

\begin{proposition} {\rm \cite[Proposition 7.8.2]{Pr-Se}} There exists a
function $\tilde{\mu}\co  Gr_{\infty}(H) \to \t\oplus \R$ such that for
any $\xi \in \t\oplus \R$, we have
$$d\langle \tilde{\mu}, \xi \rangle = \omega (X_{\xi}, \cdot)$$ where
$X_{\xi}$ denotes the infinitesimal vector field on $Gr_{\infty}(H)$
induced by $\xi$. 
\end{proposition}

The space $Gr_0(H)$ is also $T\times S^1$ invariant. Moreover, it is 
invariant under the natural action of the complexified torus $T^{\bC}$ on
$Gr(H)$, which is induced by the action of $T^{\bC}$ on $H$.  
Combined
with the action of $(S^1)^{\bC}= \bC^*$ defined in \cite[\S
7.6]{Pr-Se}, this gives an action of the complex torus $T^{\bC}\times
\bC^*$ on $Gr_0(H)$, which extends the action of $T\times S^1$.  Note that for a fixed $n$, the
identification $$W \mapsto W/ e^{in\theta}H_+ \subset
e^{-in\theta}H_+/e^{in\theta}H_+=\bC^{2nN}$$ makes the space of all
$W$ that satisfy equation~\eqref{H} into a finite-dimensional
Grassmannian, denoted
\begin{equation}\label{cal}
{\mathcal G}_n := Gr(nN, \bC^{2nN}).
\end{equation}
Thus we may think of ${\mathcal G}_n$ as a subset of
$Gr_{\infty}(H)$. By the definition of the action, it is
straightforward to see that each ${\mathcal G}_n \subseteq
Gr_{\infty}(H)$ is $T^{\bC} \times \bC^*$--invariant.  The following
proposition summarizes the facts we need about the ``Grassmannian
model" \cite[\S 8.3]{Pr-Se}, \cite[Lemma 2.4]{Liu}. 

\begin{proposition}\label{imath} 
\begin{enumerate}
\item[(a)] Let $\Omega_{\infty}(G)$ denote the space of all smooth loops $S^1\to G$.
 The map $\imath_1\co  \Omega_{\infty}(G) \to Gr_{\infty}(H)$
given by   $\imath_1(\gamma) := \gamma H_{+} $, is injective.   
We have $\mu|_{\Omega_{\infty}(G)} =\tilde{\mu} \circ
\imath_1$.

\item[(b)] We have $\Omega_{\a}(G) \subset \Omega_{\infty}(G)$ and $\imath_1 ( \Omega_{\a}(G)) \subset
Gr_{0}(H)$. Consequently,\break $\imath_1 ( \Omega_{\a}(G)) $ is contained
in $\bigcup_{n\ge 0} X_n$, where $X_n$ denotes the space of all closed
linear subspaces $W\subset H$ with the properties $$ e^{in\theta}H_+
\subset W \subset e^{-in\theta}H_+, \quad e^{i\theta}W\subset W.$$

\item[(c)] The map $W\mapsto W/ e^{in\theta}H_+$ identifies
  $X_n$ with a subvariety of the Grassmannian ${\mathcal G}_n$ of
  all complex vector subspaces of
  $e^{-in\theta}H_+/e^{in\theta}H_+ \cong \bC^{2nN}$ of dimension $nN$. The subspace
  ${\mathcal G}_n \hookrightarrow Gr_{\infty}(H)$ is $T\times S^1$
  invariant.  The restriction of the K\"ahler form on $Gr(H)$
  to ${\mathcal G}_n $ is equal
  to the K\"ahler structure induced on ${\mathcal G}_n \cong
  Gr(nN,\bC^{2nN})$ via the Pl\"ucker embedding into $\P^M$, where \(M =\)
  \(2nN\choose nN\).

\end{enumerate}
\end{proposition}

\subsection{Morse theory for the components of $\mu$}\label{subsecMorse}

In this subsection we discuss Morse theory for the components of the
moment map $\mu\co  \Omega_1(G) \to \t\oplus \R$.  
We  recall that an
element of $\t$ is called {\it regular} if it is not contained in any
of the hyperplanes $\ker \alpha \subseteq \t$, where $\alpha \in \t^*$
is a root of $G$.\footnote{We apologize to the reader for the
confusing terminology; a ``regular element of $\t$" means something
different from a ``regular value of the moment map." Unfortunately both
terms are fairly standard in their respective contexts, so we use both
here. We hope the context makes clear which we mean.}  The connected
components of $\t \setminus \bigcup_{\alpha \ {\rm root}} \ker\alpha$ 
 are called Weyl chambers.
We will need the following  result (for a proof, one can see
for instance \cite[Chapter V, Proposition 2.3]{Br-tD}).
\begin{lemma} If $X$ is an element of  $\t$ then the following two assertions are
equivalent.
\begin{itemize}
\item[(i)] The centralizer of $\exp(X)$ reduces to $T$.
\item[(ii)] The vector $X$ does not belong to any of the affine hyperplanes
$$H_{\alpha,k}:=\{ H\in \t \ | \ \alpha(H)= k\},$$
where $\alpha$ is a root and $k\in \Z$.
\end{itemize}
\end{lemma}
If $X\in \t$ is  an arbitrary regular vector, then there exists
a positive  integer  $q(X)$ such that for any 
$q \ge q(X)$, the vector $X/q$ satisfies condition (ii) from
above.   

\begin{definition}\label{definition} An element $\rho\in \t$ is called {\rm admissible}
if it is of the type $X/q$, where $X$ is a regular element which 
belongs to the integer lattice and
$q\ge q(X)$.
\end{definition}

We define a pairing on the Lie algebra $\t \oplus \R$ by taking
 the restriction of the pairing on $\g$ on the first factor, taking
 the standard inner product on the second factor, and declaring the
 two factors to be orthogonal. We denote this pairing also by $\langle
 \cdot, \cdot \rangle$. For the purposes of the Morse theoretical
 arguments in the following sections, we are
 interested in the functions $\langle \mu(\cdot), (0,1)\rangle$ and
 $\langle \mu(\cdot), (\rho,1)\rangle$, where $\rho\in \t$ is an admissible
 element. The first function is the energy functional $E$,
 namely the $S^1$--component of the moment map. The second one is, up to
 an additive constant $c=c(\rho)$, the ``tilted" energy
 $\tilde{\mathcal E}$ \cite[\S 8.9]{Pr-Se}, ie,
\begin{equation}\label{tilted}\langle \mu(\gamma), (\rho,1)\rangle = c+ \frac{1}{4\pi}
\int_0^{2\pi} \| \gamma(\theta)^{-1}\gamma'(\theta)+\rho\|^2
d\theta,\end{equation} $\gamma \in \Omega_1(G)$.  Both of these functions have a
good corresponding Morse theory with respect to the K\"ahler metric on
$\Omega_1(G)$; in particular, the downward gradient flow with respect
to both of these functions is well-defined for all time $t$, and both
the critical sets and their unstable manifolds have explicit
descriptions. The following two results are proved in \cite[\S
8.9]{Pr-Se} (see also \cite[\S 3.1]{Ko}). We denote by $\nabla f$ the
gradient vector field corresponding to a function $f$.

\begin{proposition}\label{critical} 

\begin{enumerate}
\item[(a)] Let $\rho$ be admissible vector in $\t$.
The critical
set of $\mu^{(\rho,1)}:=\langle \mu(\cdot), (\rho,1)\rangle$ is the
lattice $\check{T}$ of group homomorphisms $\lambda \co  S^1 \to T$. If
$\gamma \in \Omega_{\a}(G)$, then the solution $\phi_t(\gamma)$ of the
initial value problem
$$\frac{d}{dt}\phi_t(\gamma) = -\nabla(\mu^{(\rho,1)})_{\phi_t(\gamma)},
\quad \phi_0(\gamma) = \gamma$$ is defined for all $t\in
\R$. Moreover, $\phi_t(\gamma)\in \Omega_{\a}(G)$ for all $t\in \R$,
the limit $\lim_{t\to -\infty}\phi_t(\gamma)$ exists, and it is a
critical point.

\item[(b)] If $\lambda$ is a critical point of $\mu^{({\rho},1)}$, then the unstable manifold 
$$\C(\rho):=\{\gamma \in \Omega_{\a}\ : \
  \lim_{t\to -\infty}\phi_t(\gamma) = \lambda\}$$
  is homeomorphic to $\bC^{m(\lambda,\rho)}$, for a certain number $m(\lambda,\rho)$.

\item[(c)] We have the cell decomposition
$$\Omega_{\a}(G) = \bigcup_{\lambda\in \check{T}}\C(\rho).$$

\item[(d)] For any $\rho$ in the positive Weyl chamber,
the $\C(\rho)$ are the Bruhat cells $\C$ (as in \cite[\S 8.6]{Pr-Se}).
\end{enumerate}

\end{proposition}

We have a similar result for the other component of $\mu$, namely
$E$. The main difference is that the critical sets are no longer
isolated, but instead are diffeomorphic to coadjoint orbits of $G$. 

\begin{proposition}\label{energy} 
\begin{enumerate}
\item[(a)] The critical set of the energy functional $E =\langle
\mu(\cdot), (0,1)\rangle$ is the union of the spaces
$$\Lambda_{|\lambda|} = \{g \lambda g^{-1} \ : \ g\in G\},$$ where $
|\lambda| \in \check{T}/W$.  If $\gamma \in \Omega_{\a}(G)$, then the
solution $\psi_t(\gamma)$ of the initial value problem
$$\frac{d}{dt}\psi_t(\gamma) = -\nabla E_{\psi_t(\gamma)}, \quad \psi_0(\gamma) = \gamma$$
is defined for all $t\in \R$. Moreover, $\psi_t(\gamma)\in \Omega_{\a}(G)$ for all $t\in \R$, the limit $\lim_{t\to -\infty}\psi_t(\gamma)$ exists, and it is a critical point.

\item[(b)] For any $\lambda \in \check{T}$, the unstable manifold
  $$C_{|\lambda|}:=\{\gamma \in \Omega_{\a}\ : \ \lim_{t\to
  -\infty}\phi_t(\gamma) \in \Lambda_{|\lambda|} \}$$ is a closed
  submanifold of $\Omega_1(G)$.

\item[(c)] We have the  decomposition
$$\Omega_{\a}(G) = \bigcup_{\lambda\in \check{T}/W}C_{|\lambda|}.$$

\end{enumerate}
\end{proposition}

It is known that  many results of Morse theory
remain true for real functions on infinite dimensional Riemannian manifolds,
provided that they satisfy the condition (C) of Palais and Smale.
The latter condition is as follows (see
for instance  \cite[Chapter 9]{Pa-Te}).

\begin{definition}[Condition (C)]\label{palaissmale}
We say that a real function $f$ on a Riemannian manifold $M$ satisfies the
condition (C) if any sequence $\{s(n)\}$ of points in $M$ for which $\{(f(s(n))\}$ is
bounded and $\lim_{n\to \infty}\|\nabla (f)_{s(n)}\|=0$ has a convergent subsequence. 
\end{definition}

It turns out that the components of the moment map on
$\Omega_1(G)$ do satisfy this condition, with respect to 
the $H^1$ metric (see \fullref{subsecalgebraicloops}).
Indeed, let us consider the ``gauge" action of $L_1(G)=H^1(S^1,G)$ on $H^0(S^1,\g)$,
given by
$$\gamma \star u =\gamma u \gamma^{-1} - \gamma'\gamma^{-1},$$
for all $\gamma \in L_1(G)$ and all $u\in H^0(S^1,\g)$. This action is proper and 
Fredholm, which implies that if we regard any of its orbits as a sumbanifold of 
(the Hilbert space) $H^0(S^1,\g)$, then for any $a\in H^0(S^1,g)$ the distance function 
$f_{a}(x)=\|x+a\|^2$,  with $x$ on the orbit,
satisfy the condition (C) (see \cite[Proposition 2.16 and \S 4]{Te1}). 
The stabilizer of the constant loop $0\in H^0(S^1,\g)$ consists of all constant loops
in $G$, thus the orbit of $0$ is $L_1(G)/G=\Omega_1(G)$. It is not difficult to see that any element of the orbit can be written uniquely as $\gamma^{-1}\gamma'$,
with $\gamma\in \Omega_1(G)$. It turns out (see again \cite[\S 4]{Te1}) that the  metric on $\Omega_1(G)$ induced by its embedding in $H^0(S^1,\g)$ is just the $H^1$ metric. We summarize as 
follows.

\begin{proposition}\label{propTerng}
Regard $\Omega_1(G)$ as a Riemannian manifold with respect to the $H^1$ metric.
The map $$F\co  \Omega_1(G) \to H^0(S^1, \g),\quad F(\gamma) =
\gamma^{-1}\gamma'$$ is an isometric embedding.   For any $\rho \in \t$, the function $f_{\rho} \co  F(\Omega_1(G)) \to \R$
given by $f_{\rho}(x) = \|x+\rho\|^2$ satisfies the condition
(C). Consequently, the function $\langle \mu, (\rho,1)\rangle = c +
\frac{1}{4\pi} f_{\rho}\circ F$ (see equation \eqref{tilted}) defined on $\Omega_1(G)$ satisfies the condition (C) as well.
\end{proposition}

The fact that the components of the moment map $\mu$ on $\Omega_1(G)$
satisfy condition (C) will be the key element in our arguments in
Sections~\ref{secregular} and~\ref{secsingular}. 

\subsection{Geometric invariant theory}\label{GIT} 

In our proof of the connectivity result for the algebraic loops
$\Omega_{\a}(G)$, we will heavily use 
 the fact that in certain cases, the symplectic
quotient is the same as the geometric invariant theory quotient (in
the sense of algebraic geometry). 
The main reference for this section is the book \cite{MFK} (see also
\cite{Ki1,Ki2,Do}).  Assume that a complex torus $T^{\bC}=(\bC^*)^k$
acts analytically on a complex projective irreducible (possibly
singular) variety $X$. Let $L$ be an ample $T^{\bC}$--line bundle on $X$. A
point $x\in X$ is called {\em $L$--semistable} if there exists $m\ge 1$
and a section $s$ of the tensor bundle $L^{\otimes m}$ with $s(x)\neq
0$. The set $X^{ss}(L)$ of all semistable points in $X$ is a Zariski
open, hence irreducible, subvariety of $X$. If $Z\subset X$ is a
subvariety, then
$$Z^{ss}(L) = Z\cap X^{ss}(L).$$
Now suppose that $Y\subset \P^M$ is a smooth projective variety,
invariant under the linear action of $T^{\bC}=(\bC^*)^k$ on $\P^M$. Take
$\mu\co  Y \to \t:={\rm Lie}(T)$ to be the restriction to $Y$ of any $T$--moment
map on $\P^M$. We say that $x\in Y$ is {\em $\mu$--semistable} if
$\overline{T^{\bC}.x}\cap \mu^{-1}(0) \neq \phi$. We denote by
$Y^{ss}(\mu)$ the set of all $\mu$--semistable points in $Y$.  The
following result is Theorem 8.10 in \cite{Ki1}:

\begin{theorem}\label{kirwan} Let $Y^{\rm min}$ denote the minimal Morse--Kirwan stratum of $Y$ with respect to the function \(f = \|\mu\|^2,\) corresponding to the critical set \(f^{-1}(0) = \mu^{-1}(0).\) Then \(Y^{ss}(\mu) = Y^{\rm min}.\) 
\end{theorem}

In order to relate the $L$--semistable points and $\mu$--semistable
points, we will also need the following result of Heinzner and
Migliorini \cite{He-Mi}.
\begin{theorem}{\rm\cite{He-Mi}}\label{heinzner}\qua  There exists an ample $T^{\bC}$-line bundle $L$ such that
$Y^{ss}(L)=Y^{ss}(\mu)$.
\end{theorem}

\section{The levels of moment maps on algebraic loops}\label{secalgebraic}

The goal of this section is to prove that all non-empty levels of both
$E\co \Omega_{\a}(G) \to \R$ and $\mu\co  \Omega_{\a}(G) \to \t \oplus \R$
are connected. The essential idea of the proof is to use the
decomposition of $\Omega_{\a}(G)$ into Bruhat cells, given in
\fullref{critical}. We will first prove that the level sets on
the closure of each Bruhat cell $C_{\lambda}$ are connected, using the
fact that each such closure \(\overline{C_{\lambda}}\) can be
interpreted as a $T$--invariant subvariety of an appropriate projective
space. This allows us to use geometric invariant theory results
(Theorems~\ref{kirwan} and~\ref{heinzner}) to prove the
connectivity result on each $\overline{C_{\lambda}}$. We then use the closure
relations on the Bruhat cells to prove that the connectivity result
for each cell is sufficient to show that the level set on the union
\(\Omega_{\a}(G) = \cup_{\lambda} C_{\lambda}\) is also connected.

In order to use the results outlined in \fullref{GIT}, we must
first prove that the spaces under consideration are indeed
$T$--invariant Zariski-closed subvarieties of an appropriate projective
space. Recall that $\Omega_{\a}(G)$ contains subvarieties of finite-dimensional
Grassmannians ${\mathcal G}_n$, as explained in
\fullref{subsecloopgroups}. On the other hand, the Grassmannians
${\mathcal G}_n = Gr(n,\bC^{2n})$ may also be seen as a subvariety of a
projective space via the Pl\"ucker embedding. We first observe that
this Pl\"ucker embedding is compatible with the torus action
on $\Omega_{\a}(G)$. 

\begin{lemma}\label{lemmaPlucker}
 Let $\varphi\co  {\mathcal G}_n \hookrightarrow \P^M$ be the Pl\"ucker embedding
of the Grassmannian ${\mathcal G}_n$ where $M:=$ ${2nN}\choose{nN}$.  Then
the action of $T^{\bC}\times \bC^*$ on ${\mathcal G}_n$ extends to a linear action on $\P^M$.
\end{lemma}

\begin{proof} By definition, ${\mathcal G}_n$ is the space of all $nN$--dimensional subspaces
of
\begin{equation}\label{basis}e^{-in\theta}H_+/e^{in\theta}H_+\simeq
\langle e^{-ik\theta}b_j\ : \ -n\le k\le n-1, 1\le j \le N\rangle
\simeq \bC^{2kn},\end{equation} where $b_1, \ldots, b_N$ is the
canonical basis of $\bC^N$.  Recall that in order to define the
Pl\"ucker embedding, we specify an element of the ${\mathcal G}_n$ by
a \(2nN \times nN\) matrix $A$ of rank $nN$; such a matrix specifies
a subspace of $\bC^{2N}$ by taking the span of its columns. The
Pl\"ucker embedding then maps $A$ to $\P^M$ by taking the coordinates
in $\P^M$ to be the \(nN \times nN\) minors of $A$. This is
well-defined.

We wish to show that this standard Pl\"ucker embedding is
$T^{\bC}\times \bC^*$--equivariant, and in fact the action extends to a linear
$T^{\bC}\times \bC^*$--action on $\P^M$. By our assumption on $T$ (see
\fullref{subsecloopgroups}), $T^{\bC}$ acts on $\bC^N$ as left
multiplication by diagonal matrices. The corresponding action of
$T^{\bC}$ on the matrix $A$, which represents an element in ${\mathcal
G}_n$, is therefore also given by left multiplication by a diagonal
matrix. It is then straightforward to check that the action of
$T^{\bC}$ extends to a linear (indeed, diagonal) action on $\P^M$.

It remains to show that the extra $\bC^*$--action also extends linearly to $\P^M$. 
Since this action is to ``rotate the loop," an element $z$ of $\bC^*$ acts as follows:
$$z.(e^{ik\theta}b_j) = (ze^{i\theta} )^k b_j=z^k(e^{ik\theta}b_j),$$
$-n\le k\le n-1, 1\le j \le N$. Again, regarded as a linear transformation of $\bC^{2nN}$,
$z$ is  diagonal. We use the same argument as before to complete the proof.
\end{proof}

We are now prepared to prove the connectivity for each Bruhat cell. 
Let $\lambda \in \check{T}$. There exists $n\ge 1$ such that
$\imath_1(\overline{\C})\subset X_n$ (see eg, \cite{At-Pr}). Here
$\overline{\C}$ represents the closure of $\C$ in $\Omega_{\a}(G)$ in
the {\em direct limit} topology, as explained in the Introduction. 
The $T\times S^1$--equivariance of the Pl\"ucker
embedding in the previous lemma allows us to analyze the level set of
the moment map $\mu$ on $\overline{\C}$ as that of a moment map
defined on a finite-dimensional space, namely the Grassmannian. The
next result is a direct consequence of \fullref{imath} (see
also \cite[Lemma 2.4]{Liu}).

\begin{lemma}\label{diagram}
The following diagram is commutative.
\[
\xymatrix@!=1pc{
& \Omega_{\a}(G) \ar[dr]^{\mu} & \\
\overline{\C} \ar[dr]_{\imath_1} \ar[ur]^{\imath} & & \t \oplus
\R \\
& {\mathcal G}_n \ar[ur]_{\mu_n} & \\
}
\]
Here $\imath\co \overline{\C} \to \Omega_{\a}(G)$ denotes the  inclusion map, and 
$\mu_n \co  {\mathcal G}_n \to \t \oplus \R$ is the moment map corresponding
to the action of $T\times S^1$ on ${\mathcal G}_n$ and the symplectic form on ${\mathcal G}_n$ 
induced by the Pl\"ucker
embedding $\varphi\co  {\mathcal G}_n \hookrightarrow \P^M$. 
\end{lemma}

\begin{proof} By \fullref{imath} (a), we have
$$\mu\circ \imath = \mu|_{\overline{\C}}=\tilde{\mu}\circ \imath_1 |_{\overline{\C}}.$$
We only need to show that $\tilde{\mu}|_{{\mathcal G}_n} = \mu_n$.
We have already seen in \fullref{imath} that the symplectic
structures on ${\mathcal G}_n$ coming from the inclusions \({\mathcal
  G}_n \into \P^M\) and \({\mathcal G}_n \into Gr(H)\) 
   are
equal. Hence  it is sufficient to
see that the two inclusions are \(T\times S^1\)--equivariant. For the
first inclusion, this statement is the content of \fullref{lemmaPlucker}, and
the latter was observed in \fullref{subsecloopgroups} (see \fullref{imath} (c)). 
\end{proof}

We can now prove the connectivity result on each closed Bruhat cell, using
the geometric invariant results of \fullref{GIT}. 

\begin{lemma}\label{path} 
If \(\lambda \in \check{T},\) then for any value \(a \in \t \oplus
\R,\) the level set \(\mu^{-1}(a) \cap \overline{\C}\) is either empty
or path
connected.
\end{lemma}

\begin{proof} 
By \fullref{diagram}, we have $$\mu^{-1}(a) \cap
\overline{\C}=\mu_n^{-1}(a)\cap \overline{\C}.$$ Hence it suffices to
prove that the RHS is connected.  We consider the function $f\co 
{\mathcal G}_n \to \R$, $f(x) = \|\mu_n(x)-a\|^2$, as well as its
gradient flow $\Phi_t$. We will prove first that $\Phi_t$ leaves
$\overline{\C}$ invariant.  It suffices to show that
the gradient flow of $f$ leaves any $T^{\bC}\times \bC^*$ orbit in
${\mathcal G}_n$ invariant, since $C_{\lambda}$ is $T^{\bC}\times \bC^*$--invariant.
Recall that the gradient vector field of $f$ is given by 
$${\rm grad}(f)_x = J (\mu_n(x)-a).x,$$ where the dot denotes the
infinitesimal action of ${\rm Lie}(T\times S^1)=\t
\oplus \R$ on $\mathcal{G}_n$ at the point \(x \in \mathcal{G}_n\), 
and $J$ denotes the complex structure. 
We deduce that the vector field
${\rm grad}(f)$ is tangent to a $T^{\bC}\times \bC^*$ orbit,
hence any integral curve of ${\rm grad}(f)$ remains in an orbit. 

Recall that the function $f$ induces a stratification of ${\mathcal
G}_n$ (\cite{Ki1}) as follows. Let \(\{C_j\}_{j=0}^{\ell}\) be the 
critical
sets of $f$ in ${\mathcal G}_n$, ordered by their values $f(C_j)$, 
where $C_0=f^{-1}(0)=\mu_n^{-1}(a)$ is the
minimum.  To each $C_j$ corresponds the stratum
$$S_j=\{x\in {\mathcal G}_n\ : \ \lim_{t\to \infty} \Phi_t(x) \in
C_j\}.$$ 
Each $S_j$ is a complex
submanifold of ${\mathcal G}_n$, which contains $C_j$ as a deformation
retract \cite{Wo,Le,MFK}. 
The unique open stratum $S_0$ coincides with the set
$({\mathcal G}_n)^{ss}(\mu_n-a)$ of all semi-stable points in
${\mathcal G}_n$ (\fullref{kirwan}).  By \fullref{heinzner}, there exists an ample $T^{\bC}\times \bC^*$--line
bundle $L$ such that $$({\mathcal G}_n)^{ss}(\mu_n-a) = ({\mathcal
G}_n)^{ss}(L).$$ Since $\overline{\C}$ is a subvariety of ${\mathcal
G}_n$, the set of all semi-stable points of $\overline{\C}$ is
$$(\overline{\C})^{ss}(L) = ({\mathcal G}_n)^{ss}(L)\cap \overline{\C}
 =S_0\cap \overline{\C}.$$ Since $\overline{\C}$ is irreducible (see
 eg, \cite{Mi}), $(\overline{\C})^{ss}(L)$ is irreducible as well. Since
 any irreducible complex analytic projective variety has a resolution
 of singularities, $(\overline{\C})^{ss}(L)$ is path-connected in the
 usual differential topology of $\P^M$. 
 We 
 deduce that $S_0\cap \overline{\C}$ is path-connected. Because
 $\overline{\C}$ is invariant under $\Phi_t$, the space
 $S_0\cap\overline{\C}$ contains $C_0\cap \overline{\C}=\mu_n^{-1}(a)
 \cap \overline{\C}$ as a
 deformation retract. Thus the latter is path-connected as well.
\end{proof}

We now use
\fullref{path} and the closure relations for the Bruhat cells
to prove the connectivity result for all of $\Omega_{\a}(G)$. 

\begin{proposition}\label{muconnected}
For any $a\in \t\oplus \R$, the set $\mu^{-1}(a) \cap \Omega_{\a}(G)$ is
empty or connected.
\end{proposition}

\begin{proof} 
Let $\gamma_1$ and $\gamma_2$ be in $\mu^{-1}(a) \cap
\Omega_{\a}(G)$. Since the union of the Bruhat cells equals
$\Omega_{\a}(G)$, there exist $\lambda_1,\lambda_2$ in the integer lattice $\check{T}$ such that
\(\gamma_1\in C_{\lambda_1}, \gamma_2 \in C_{\lambda_2}.\) 
Let us consider the  ordering $\le$ on $\check{T}$ given by $\lambda \le \nu$ if and only if $C_{\lambda}\subset \overline{C_{\nu}}$. One can show that this is the same as the Bruhat ordering on $\check{T}$ (see for instance \cite[Theorem 1.3]{Mi}). This is 
induced (see \cite[\S 1.3]{Ku}) by the Bruhat ordering on the affine Weyl group 
$W_{{\rm aff}} = \check{T}\rtimes W$, which is a Coxeter group, via the obvious identification $\check{T}=W_{{\rm aff}}/W$. One of the properties of the Bruhat ordering on Coxeter groups is that for any two elements there exists a third one which is ``larger" than both of the previous two (see for instance \cite[Lemma 1.3.20]{Ku}). In our context this implies that there exists $\nu \in \check{T}$  such that $\lambda_1\le \nu$ and $\lambda_2\le \nu$,
that is,
\(C_{\lambda_1} \subseteq \overline{C_{\nu}}\) and \(C_{\lambda_2} \subseteq \overline{C_{\nu}}\).
 We apply
\fullref{path} and deduce that $\mu^{-1}(a) \cap \Omega_{\a}(G)$ is path
connected, hence connected.
\end{proof}

By using exactly the same methods, we can prove the following proposition.

\begin{proposition}\label{econnected} For any $a\in\R$, the energy
  level set $E^{-1}(a) \cap \Omega_{\a}(G)$ is either empty or 
  connected.
\end{proposition}

\section{The regular levels of the moment map on $\Omega(G)$}\label{secregular}

In this section we prove the connectivity result for $\Omega_1(G)$, using the
topology induced from $H^1$. By \fullref{muconnected},
we know that $\mu^{-1}(a)\cap \Omega_{\a}$ is a connected subspace of $\Omega_{\a}(G)$
in the direct limit topology.  As observed in \fullref{secintro}, the
direct limit topology is strictly finer than the topology induced from
$\Omega_1(G)$, so $\mu^{-1}(a) \cap \Omega_{\a}(G)$ is also connected
in the $H^1$ topology. For the rest of the section, we consider only
the $H^1$ topology. 

In order to prove the connectivity result for $\Omega_1(G)$, we will
heavily use  the connectivity result for $\Omega_{\a}(G)$, proven in
the previous section. We also use that
$\Omega_{\a}(G)$ is dense in $\Omega_1(G)$ in the $H^1$
topology \cite[Theorem 2]{At-Pr}.  
Since the closure of a connected subset is also connected,
in order to prove \fullref{secondmain}, it therefore suffices to
prove the following proposition. 

\begin{proposition}\label{dense} If $a \in \t \oplus \R$ is a regular value of $\mu$, then
$\mu^{-1}(a)\cap \Omega_{\a}(G)$ is dense in $\mu^{-1}(a)$ in the
  $H^1$ topology. 
\end{proposition}

\begin{remark}
The reader may notice that the above proposition has an additional
hypothesis of regularity of the value, which was not necessary for the
case of $\Omega_{\a}(G)$. We need this hypothesis to use the 
Morse-theoretic arguments below. We are not aware of other methods to
prove \fullref{secondmain}, although the question of singular
values is still of interest. We address the case of singular values
for $E$ in the next section. 
\end{remark}

We now concentrate on the proof of 
\fullref{dense}.  Before we proceed, we warn the reader that
throughout the rest of this section, the loops will be generically denoted by $x$,
rather than $\gamma$.  In order to prove the density,
we will construct for any point \(x_0 \in \mu^{-1}(a)\) a
sequence contained in \(\mu^{-1}(a) \cap \Omega_{\a}(G)\) converging
to it. We accomplish this by first using the fact that
$\Omega_{\a}(G)$ is dense in $\Omega_1(G)$ to find a sequence
\(x(r) \to x_0 \in \mu^{-1}(a),\) where \(x(r) \in \Omega_{\a}(G).\) We
then use the Morse flow with respect to independent components of the
moment map $\mu$, restricted to $\Omega_{\a}(G)$, to produce a new
sequence \(y(r) \to x_0,\) where now \(y(r) \in \mu^{-1}(a) \cap
\Omega_{\a}(G).\) It is precisely this construction -- which uses
Morse flows with respect to components of the moment map $\mu$ --
where we need the Palais--Smale condition (C). 

We will need the following. 

\begin{lemma}\label{vector} Let $a$ be a regular value of $\mu\co 
  \Omega_1(G) \to \t \oplus \R$, and let \(x \in \mu^{-1}(a).\) 

\begin{enumerate}
\item[(a)] The map $P\co  \t\oplus \R\to T_{x}(\Omega_1(G))$, $\xi \mapsto
\nabla\langle \mu(\cdot),\xi\rangle(x)$ is $\R$--linear and injective.
Denote by $V:=P(\t\oplus \R)\subset T_{x}(\Omega_1(G))$ the image of this map.

\item[(b)] Denote $k-1:=\dim (\t)$. There exists a basis $\xi_1,
  \ldots,\xi_k$ of $\t\oplus \R$ such that 
$\xi_1=(0,1)$, $\xi_j=(\rho_j, 1)$, where $\rho_2,\ldots, \rho_k$
are admissible (in the sense of \fullref{definition}) 
  elements of $\t$ with the 
following property: let 
$v_i:= \nabla\langle \mu(\cdot),\xi_i\rangle(x)$, and $p_i \co V \to V$ the orthogonal
projection on the hyperplane orthogonal to $v_i$. Then 
$$ g(v_j , p_1\circ \ldots \circ p_{j-1} (v_j))
\neq 0,$$
for all $j\ge 2$. Here $g$ denotes the K\"ahler metric on
$\Omega_1(G)$. 
\end{enumerate}

\end{lemma}

\begin{proof}
We first prove claim (a).  Since $a$ is a regular value, \(T\times
S^1\) acts locally freely on the level \(\mu^{-1}(a)\), hence the
linear map $\xi \mapsto \xi_{x}$ from $\t \oplus \R \to
T_x\Omega_1(G)$ is injective. The vector $\xi_{x}$ is the
Hamiltonian vector field of $\langle \mu(\cdot),\xi\rangle$ at $x$.
By properties of Hamiltonian vector fields on a K\"ahler manifold,
\(J_x \xi_x = \nabla \langle \mu(\cdot), \xi\rangle,\) where $J_x$
denotes the complex structure at $x$. We then have that \(P(\xi) = J_x
\xi_x,\) and since $J_x$ is a linear isomorphism at $T_xM$, the claim
follows.

We now prove claim (b). 
 We identify $V$ with  $\t\oplus \R$. This gives the identifications 
$$v_1:=P(\xi_1) =\xi_1:=(0,1),\quad \mbox{and} \quad
 v_j:=P(\xi_j)=\xi_j:=(\rho_j,1),
\quad 2 \le j \le k. $$
The space $\t\oplus \R$  inherits from $T_{x_0}(\Omega(G))$ an
inner product via the map $P$. By abuse of notation we denote this
metric on $\t \oplus \R$ also by $g$. The argument is more easily
 given on $\t \oplus \R$, but since the metric $g$ is the pullback
 from $T_x\Omega_1(G)$, the claim follows from the following
 construction of the $\rho_i$. 

We inductively construct  admissible vectors $\rho_2, \ldots, \rho_k\in
\t$ with the property that for any $\ell \ge 2$ we have
\begin{itemize}
\item[(i)] the vectors $(0,1), (\rho_2,1), \ldots, (\rho_{\ell},1)$  are linearly independent, 
\item[(ii)]  $g(v_{\ell} , p_1\circ \cdots \circ p_{\ell-1} (v_{\ell}))
\neq 0,$
\item[(iii)] $\label{gammath}g( (0,1), (\rho_{\ell},1))\neq 0.$
\end{itemize}
Let \(2 \leq \ell \leq k-1.\) If \(\ell \geq 3\) we assume that \(\{\rho_2,
 \ldots, \rho_{\ell-1}\}\) have already been constructed. 
Now consider the  function
 $f\co \t\to \R$,  given by
\begin{align} f(\rho)=& g( ( \rho,1), p_1 \circ \cdots \circ p_{\ell-1}(\rho,1 ))\nonumber \\{}
=& g((\rho,0), p_1 \circ \cdots \circ p_{\ell-1}(\rho,1)) \nonumber \\ =&  
g(( \rho,0), p_1 \circ \cdots \circ p_{\ell-1}(0,1)) + g((\rho,0),
p_1 \circ \cdots \circ p_{\ell-1}(\rho,0)).\nonumber
\end{align}
The two components of $f$ described by the previous equation are homogeneous polynomial
functions of degree 1, respectively 2, on $\t$. 

We will prove that $f$ is not identically zero. To this end, it is sufficient to prove that
the degree 1 component of $f$ 
is not identically zero. More specifically, we will prove the following claim.

{\bf Claim}\qua The function 
$\rho \mapsto g( (\rho,0), p_1\circ\cdots \circ
p_{\ell-1}(0,1))$, $\rho \in \t$
is not identically zero.

{\bf Proof of the claim}\qua First we  show that $p_1 \circ \ldots \circ
p_{\ell-1}(0,1) 
=p_1 \circ \ldots \circ p_{\ell-1}(\xi_1)$ 
is different from the $0$ vector.
To do that, we take into account that
\begin{equation}\label{pj} p_j(\xi) = \xi-\frac{g(\xi, \xi_{j}) }{g (\xi_{j}, \xi_{j})}\xi_j \end{equation}
where $1\le j \le \ell -1$ and $\xi\in \t\oplus \R$.  From here, a
straightforward computation that involves applying a projection of
the form \eqref{pj} successively $\ell -2$ times shows that
$p_2\ldots p_{\ell-1}(\xi_1)$ is a linear combination of $\xi_1,\ldots
,\xi_{\ell-1}$, where the coefficient of $\xi_{\ell-1}$ is
$\frac{g(\xi_1, \xi_{\ell-1}) }{g (\xi_{\ell-1}, \xi_{\ell-1})}$. The
latter number is nonzero, by the induction hypothesis, hence
$p_2\ldots p_{\ell-1}(\xi_1)$ cannot be collinear to $\xi_1$. This
implies that $p_1\ldots p_{\ell-1}(\xi_1)\neq 0$, as stated above. Now
if the function mentioned in the claim was identically 0, from the
fact that
$$g((0,1), p_1\ldots p_{\ell-1}(0,1))=0$$ and the
nondegeneracy of the metric $g$ on $\t\oplus \R$, we  deduce $p_1\ldots p_{\ell-1}(0,1)=0$,
which is false.  The claim is now proved.

We now note that 
condition (ii) 
is equivalent to the statement $f(\rho_{\ell})\neq 0$.
Similarly, conditions (i)
and (iii) can be phrased as the non-vanishing of certain non-zero
polynomials $f_1$ and $f_3$ on $\t$. 
In particular, condition
(i) is equivalent to the non-vanishing of all the $\ell\times \ell$
minors (and hence of their product, which we call $f_1$) of the matrix formed by the
coordinates of the vectors $ (0,1), (\rho_2,1), \ldots,
(\rho_{\ell},1)$. Similarly, condition (iii) is equivalent to 
the non-vanishing of the degree 1 polynomial $f_3(\rho):= g((0,1), (\rho,1))$.  
Hence conditions (i), (ii), and (iii) are satisfied by any
$\rho_{\ell}$ for which 
$f_0(\rho_{\ell +1}) \neq 0$, where we define $$f_0(\rho):=
f_1(\rho)f(\rho)f_3(\rho).$$ Since $f$ was shown above to be not
identically zero, and $f_1, f_3$ are also both not identically zero,
$f_0$ is also not identically zero. 

Let $\Lambda$ denote the integer lattice in $\t$. 
Since a non-zero polynomial cannot vanish on all points of a full 
lattice, there exists $X\in \Lambda$ regular 
such that $f_0(X)\neq 0$. 
Moreover, there exists an integer $q\ge q(X)$  (see subsection 2.2 for the
definition of $q(X)$) such that $f_0(X/q) \ne 0$. This follows from
the fact that the polynomial in one variable $t$ given by $p(t) :=
f_0(tX)$ also cannot be identically zero, since $p(1)=f_0(X) \neq 0$. 
We set 
$\rho_{\ell}:=X/q$, and the lemma is proved.
\end{proof}

\begin{remark}
Notice that the first vector $\xi_1$ in this basis is deliberately
chosen so that the corresponding component \(\mu^{\xi_1}\) is exactly
the energy functional $E$. This is not necessary for the argument, but
makes it evident how to apply a similar argument in \fullref{secsingular} for
just the energy functional. 
\end{remark}

Denote by $h_1 =\langle \mu(\cdot),\xi_1\rangle , \ldots, h_k=\langle
\mu(\cdot),\xi_k\rangle$ the components of the moment map 
$\mu$ corresponding to the
$\xi_k$. By \fullref{propTerng}, each $h_j$ satisfies
condition (C) of Palais--Smale with respect to the $H^1$ metric on $\Omega_1(G)$. 
Also let $a_1 =\langle
a,\xi_1\rangle , \ldots, a_k=\langle a,\xi_k\rangle$ be the coordinates
of $a$ with respect to the basis $\{\xi_k\}$. Note that
$$\mu^{-1}(a) = h_1^{-1}(a_1)\cap \ldots \cap h_k^{-1}(a_k).$$
Let $x_0$ be a point in  $\mu^{-1}(a)$, and fix $j\in \{1,2,\ldots, k\}$. 
We denote by $g( \ , \ )$ the K\"ahler metric and by
$\langle \ , \rangle _1$ the $H^1$ Riemannian 
metric on $\Omega_1(G)$. 
We also denote by $\nabla h_j$ the gradient vector field
of $h_j$ with respect to the K\"ahler metric. Since $x_0$ is a regular
value of $h_j$, we have that \(\nabla h_j(x_0) \neq 0,\) so 
by a continuity argument
there exists an open neighbourhood $U$
of $x_0$ in $\Omega_1(G)$ and a positive number $M$ such that the
vector field on $U$ defined as 
$$Y_j:=-\frac{1}{g(\nabla h_j, \nabla h_j)}\nabla h_j $$
satisfies $\|Y_j(x)\|_1\le M$ for any $x\in U$. We now show that $Y_j$ can be extended 
to a vector field on all of $\Omega_1(G)$, with bounded $H^1$--length. Consider
a coordinate system around $x_0$ which maps $x_0$ to $0$. We may
assume without loss of generality that 
the field $Y_j$  has bounded $H^1$--length on the ball $B(0,r)$ in this
coordinate system. 
We consider a smooth function 
$g$ on the coordinate system which is equal to 1 on $B(0,r/3)$ and $0$ outside the ball
$B(0,r/2)$. 
Define the vector field $Y'_j := g Y_j$ on this coordinate
system. By extending further by $0$ to all of $\Omega_1(G)$, we may
take $Y'_j$ to be defined on all of $\Omega_1(G)$. 
We will need the following result. 

\begin{lemma} There exists an open neighbourhood $U_0$ of $x_0$ and a number
$\varepsilon >0$ with the property that for each $j=1,2,\ldots, k$, there exists  a one-parameter
group of automorphisms $\Phi^j_t$, $t\in \R$, of $\Omega_1(G)$, such
that for any $t\in (-\varepsilon,\varepsilon)$ and for any $x\in
 U_0$, we have
\begin{equation} h_j(\Phi^j_t(x))=h_j(x) - t.\end{equation}
The group $\Phi^j_t$ leaves $\Omega_{\a}(G)$ invariant.
\end{lemma}    
\begin{proof} 
Because the metric $H^1$ is complete and the vector field $Y'_j$ on
$\Omega_1(G)$ defined above has bounded $H^1$--length, we deduce that
$Y'_j$ is completely integrable (see \cite[Corollary 9.1.5]{Pa-Te}).  Let
$\Phi_t^j$, $t\in \R$, be the flow given by $Y'_j$.  By continuity, 
there exists $\varepsilon >0$ and $U_0\subset U$ such that
for any $t\in (-\varepsilon, \varepsilon)$ and any $x\in U_0$, we have
$$\Phi_t^j(x)\in U.$$
Consider the function $$t\mapsto \gamma_j(t):=h_j(\Phi^j_t(x)).$$
For  $t\in (-\varepsilon, \varepsilon)$ and  $x\in U_0$ we have
$$\frac{d\gamma_j}{dt}=
-g( \nabla h_j(\Phi^j_t(x)), \frac{1}{g(\nabla h_j(\Phi^j_t(x)),\nabla h_j(\Phi^j_t(x)))}\nabla h_j(\Phi^j_t(x))) =-1.$$
Since $\gamma_j(0)=h_j(x)$, we deduce that $\gamma_j(t) = h_j(x)-t.$   

We now prove  the last statement of the lemma. We will argue the
result only for $j\ge 2$, and simply note that a similar argument can be used for
$j=1$. In order to simplify the notation, we will assume that $\rho_j$
is in the positive Weyl chamber, so that $\C(\rho_j)=\C$ (see
\fullref{critical}). From the construction of
$Y'_j$ preceding the lemma, we can  see that for any $x\in \Omega_1(G)$, 
the vector $(Y'_j)_x$ 
is a multiple of $\nabla h_j(x)$. In
particular, if $\C$ is a Bruhat cell, then for any $x \in \C$, the
vector $(Y_j)_x$ is in $T_x \C$ (see \fullref{critical}). The restriction of $Y_j$
to $\C$ has bounded $H^1$ length, hence it generates a 1--parameter group of
diffeomorphisms $\tilde{\Phi}^j_t$, $t\in \R$, of $\C$ (again by \cite[Corollary 9.1.5]{Pa-Te}).
 For any $x\in \C$, the curves $\tilde{\Phi}^j_t(x)$ and
$\Phi^j_t(x)$ are integral curves of $Y'_j$ in the
Hilbert manifold $\Omega_1(G)$ with the same initial condition at
$t=0$. By \cite[Chapter IV, \S 2, Theorem 2]{La}, we deduce that
$\Phi^j_t(x)=\tilde{\Phi}^j_t(x)$ for all $t\in \R$. In other words, $\Phi^j_t$
leaves the Bruhat cell $\C$ invariant, for any $\lambda\in
\check{T}$. 
Hence it leaves $\Omega_{\a}(G)$ invariant, as desired. 
\end{proof}

Let us consider the map $\pi_j \co  U_0\cap h_j^{-1}(a_j-\varepsilon, a_j
+\varepsilon) \to h_j^{-1}(a_j)$ given by
$$\pi_j (x)= \Phi^j_{h_j(x)-a_j}(x).$$
We wish to compose these maps $\pi_j$ in order to take a sequence
in $\Omega_{\a}(G)$ to a sequence contained in $\mu^{-1}(a) \cap
\Omega_{\a}(G)$. Complications arise, however, in that there is no
guarantee that a projection \(\pi_j\) leaves \(h_i^{-1}(a_i)\)
invariant for \(i \neq j.\) In order to get around this difficulty, we
will need the following lemma. 

\begin{lemma}\label{A}
For any $1\le j \le k$, the map
$$A_j(x,t) := \pi_1\circ \ldots \circ \pi_{j-1} (\Phi^j_t(x))$$
maps an open neighbourhood of $(x_0,0)$ in $h_1^{-1}(a_1) \cap \ldots \cap h_j^{-1}(a_j)\times
(-\varepsilon, \varepsilon)$ diffeomorphically onto an open neighbourhood
of $x_0$ in $h_1^{-1}(a_1) \cap \ldots \cap h_{j-1}^{-1}(a_{j-1})$.
\end{lemma}

\begin{proof}  First of all, we note that $A_j$ is a well-defined  map
from a certain open neighbourhood of $(x_0,0)$ in $h_1^{-1}(a_1) \cap \ldots \cap h_j^{-1}(a_j)\times
(-\varepsilon, \varepsilon)$ to $h_1^{-1}(a_1) \cap \ldots \cap h_{j-1}^{-1}(a_{j-1})$.
We will show that the latter map is a local diffeomorphism at $(x_0,0)$. To this end, 
we consider the map $$B_j(x,t) := \pi_2\circ \ldots
 \circ \pi_{j-1} (\Phi^j_t(x))$$ from an open neighbourhood
of $(x_0,0)$ in $h_2^{-1}(a_2) \cap \ldots \cap h_j^{-1}(a_j)\times
(-\varepsilon, \varepsilon)$ to $h_2^{-1}(a_2) \cap \ldots \cap
 h_{j-1}^{-1}(a_{j-1})$.  We have $A_j = \pi_1\circ
 B_j|_{h_1^{-1}(a_1) \cap \ldots \cap h_j^{-1}(a_j)}$ and
 consequently the following diagram is commutative:
\[
\xymatrix @C=2.5pc {
\left(\cap_{\ell=2}^j \ker (dh_{\ell})_{x_0}\right) \oplus \R \ar[rr]^-{(dB_j)_{(x_0,0)}} &&
\cap_{\ell=2}^{j-1} \ker (dh_{\ell})_{x_0} \ar[d]^{(d\pi_1)_{x_0}} \\
\left(\cap_{\ell=1}^j \ker (dh_{\ell})_{x_0}\right) \oplus \R \ar[rr]^-{(dA_j)_{(x_0,0)}} 
\ar[u]^{\imath} && \cap_{\ell=1}^{j-1} \ker (dh_{\ell})_{x_0}
}
\]
Since $(dB_j)_{(x_0,0)} (v,0) =v$, we only need to prove that
$(dA_j)_{(x_0,0)} (0,1) \notin \ker (dh_j)_{x_0}$. This is true, because
\fullref{vector} (b) says that
\begin{equation}\label{zeroagain}g( \nabla(h_j)_{x_0} ,
(d\pi_1)_{x_0}\circ \ldots \circ (d\pi_{j-1})_{x_0} (\nabla(h_j)_{x_0} ))
\neq  0.\end{equation}
\end{proof}

We may now prove the main result of this section.

\begin{proof}[Proof of \fullref{dense}]
Let \(x_0 \in \mu^{-1}(a).\) We wish to show that \(\mu^{-1}(a) \cap
\Omega_{\a}(G)\) is dense in \(\mu^{-1}(a).\) To this end it would
suffice to construct a sequence \(y(r) \in \mu^{-1}(a) \cap
\Omega_{\a}(G)\) such that \(y(r) \to x_0\) in the $H^1$ topology. 
By \cite[Proposition 3.5.3]{Pr-Se}, we know that $\Omega_{\a}(G)$ is dense in
$\Omega_1(G)$, so 
there exists a sequence $x_0(r)$ with $x_0(r) \in \Omega_{\a}(G)$, $x_0(r)\to x_0$.
Then $$x_1(r):= \pi_1(x_0(r))$$
has the property $x_1(r)\in h_1^{-1}(a_1)\cap \Omega_{\a}(G)$ and
$x_1(r)\to \pi_1(x_0)=x_0$. Note that since \(x_0 \in h_i^{-1}(a_i) \quad \forall
i,\) the projection $\pi_i$ fixes $x_0$ for all $i$. 
Suppose now we have a sequence \(x_{\ell}(r) \to x_0,\) where
\(x_{\ell}(r) \in h_1^{-1}(a_1) \cap \cdots \cap
h_{\ell}^{-1}(a_{\ell}).\) From \fullref{A}, $A_{\ell}$ is a local
diffeomorphism, and hence by taking an appropriate subsequence, we may define
$x_{\ell+1}(r)$ by the relation 
$$x_{\ell}(r) = A_{\ell}(x_{\ell+1}(r), t_r)$$ where $x_{\ell+1}(r)
\in h_1^{-1}(a_1)\cap \cdots \cap h_{\ell+1}^{-1}(a_{\ell+1}) \cap
\Omega_{\a}(G)$. Because $A_{\ell}$ is a local diffeomorphism around
$(x_0,0)$ (see \fullref{A}), from $A_{\ell}(x_{\ell+1}(r), t_r) \to A_{\ell}(x_0,0)$ we
deduce $x_{\ell+1}(r) \to x_0$.  Continuing, we may set $y(r):=x_k(r)$
(recall that $k -1$ is the dimension of $T$). In this way \fullref{dense} is completely proved.\end{proof}

We may immediately conclude that \fullref{secondmain} holds. 

\begin{proof}[Proof of \fullref{secondmain}] 
We have just shown that $\mu^{-1}(a) \cap \Omega_{\a}(G)$ is dense in
$\mu^{-1}(a) \subseteq \Omega_1(G)$ in the $H^1$ topology. Since
$\mu^{-1}(a) \cap \Omega_{\a}(G)$ is connected in the direct limit
topology on $\Omega_{\a}(G)$, and the direct limit topology on
$\Omega_{\a}(G)$ is finer than the subset topology induced from the $H^1$
topology on $\Omega_1(G)$ (see \fullref{finer}), we may conclude that \(\mu^{-1}(a) \cap
\Omega_{\a}(G)\) is also connected in the $H^1$ topology, considered
as a subset of $\Omega_1(G)$. The closure of the connected set is also
connected, and $\mu^{-1}(a)$ is closed, so we conclude that
$\mu^{-1}(a)$ is connected, as desired. \end{proof}

By repeating the same argument using only the function $h_1 = E$, we
also immediately obtain the following, which is a special case of
\fullref{zeromain}. 

\begin{proposition}\label{propEregular}
For any  $a\in \R$ which is a regular value of $E\co \Omega_1(G) \to \R$, the
preimage $E^{-1}(a)$ is empty or connected.
\end{proposition}

\section{The singular levels of the energy on $\Omega_1(G)$}\label{secsingular}

\fullref{propEregular} of the previous section already proves
the connectivity for level sets of regular values for the energy
functional. 
Hence, in order to prove completely \fullref{zeromain}, it only remains to consider
the case of a {\em singular} level set of the energy functional $E$. 
We prove this using essentially the same ideas as in
\fullref{secregular}, except that we must remove a certain subset
of the critical points in order to force the regularity condition. It
is then also necessary to prove that the condition (C) of Palais--Smale
holds even after we have removed this subset. This is the main
technical result of this section. 

We first observe that, 
as in \fullref{secregular}, it suffices to prove the following. 

\begin{proposition}\label{denseon} Let $a \in  \R$ be a singular value
  of $E\co \Omega_1(G)\to \R$. Then 
$E^{-1}(a)\cap \Omega_{\a}(G)$ is dense in $E^{-1}(a)$.
\end{proposition}

Following the ideas of the previous section, 
we will prove below that if $x_0\in E^{-1}(a)$, then there exists a sequence
$\{x(r)\}$ with $x(r) \in E^{-1}(a)\cap\Omega_{\a}(G)$ and such that
$x(r) \to x_0$.  Denote by $K$ the (closed) subset of $\mu^{-1}(a)$
consisting of all singular points of $E$ contained in
$E^{-1}(a)$. This $K$ is of the type $\Lambda_{|\lambda|}$ in
\fullref{energy}, so it is contained in $\Omega_{\a}(G)$.
We may assume without loss of generality that $x_0\notin K \subseteq \Omega_{\a}(G)$, since
otherwise there trivially exists a sequence in $\Omega_{\a}(G)$
converging to $x_0$ (eg, the constant sequence).

We now fix, for the rest of the discussion, such a point \(x_0 \in
E^{-1}(a), x_0 \notin K.\) We must show that singular values of the energy functional are
well-behaved. We have the following two lemmas. 

\begin{lemma}\label{lcritisol}
The critical values of $E\co  \Omega_1(G) \to \R$ are isolated in $\R$.
\end{lemma}

\begin{proof} 
As shown in \cite{At-Pr}, the critical values are
$$E(\alpha)=\frac{1}{2}|\alpha|^2,$$ where $\alpha$ is in the integer
lattice $\check{T}=\ker(\exp \co  \t \to T)\subset \t$. We claim that the
set of all $E(\alpha)$ from above has no limit points. Suppose not;
then such a limit point would occur within a ball of radius $R <
\infty$.  But $\{ X \in \t: |X| \le R\} $ is compact, and therefore
only contains finitely many points in the integer lattice $\check{T}$.
\end{proof}

\begin{lemma}\label{lfi.tw} 
There exists a closed neighborhood $K_0$ of $K$ in 
$\Omega_1(G)$ such that $x_0 \notin K_0$ and 
the only critical points of $E$ in $K_0$ are those in $K$.
\end{lemma}

\begin{proof}
There exist  open neighborhoods $U$ (respectively $V$) 
of $K$  (respectively of $x_0$) in $\Omega_1(G)$ such that $U\cap V =\phi$ (because $\Omega_1(G)$ is a
normal topological space).  We choose $\varepsilon$ such that the only
critical points in $E^{-1}([a - \varepsilon, a + \varepsilon])$ are
those in $E^{-1}(a)$ (see \fullref{lcritisol}). We define $K_0$ the
closure of $U\cap E^{-1}(a - \varepsilon, a + \varepsilon) $.
\end{proof}

We consider the infinite dimensional manifold
$M:= \Omega_1(G)\setminus K_0$ (which 
is a Hilbert manifold) and the restriction of $E$ to it. 

\begin{lemma} \label{lcondc}
The function $E$ restricted to $M$ satisfies condition (C)
with respect to the restriction of the $H^1$ metric to $M$.
\end{lemma}
\begin{proof}  Condition (C) is given in \fullref{palaissmale}. 
       Let  $\{s(n)\}$  be a sequence in  $M$  with  $|E(s(n))|$ bounded
and  $\|\nabla E_{s(n)}\|_1 \to 0$. Since  $E \co  \Omega_1(G) \to \R$ does satisfy (C),
we deduce that there exists  $p \in \Omega_1(G)$  with  $s(n_k) \to  p$.
But  $p$  cannot be in  $K_0$. Indeed,  $p$  is necessarily a critical point,
so if  $p \in K_0$, then  $p \in K$ (as in the statement of 
\fullref{lfi.tw}); this
is not possible, because all  $s(n_k)$  are contained in the complement
of a neighborhood of  $K$, hence the limit must stay out of that
neighborhood. The lemma is proved.
\end{proof}

We may now prove the proposition. 

\begin{proof}[Proof of \fullref{denseon}]
For the  function $E$ restricted to $M$, $a$ 
is a {\it regular}
value. 
By using the methods from Section 4, we can construct a sequence $\{x(r)\}$ with 
$x(r) \in \Omega_{\a}(G)\cap M \cap E^{-1}(a)$ and $x(r) \to x_0$. 
The crucial fact is that $\nabla E$ is always tangent 
to Bruhat cells, hence the flow $\Phi_t$ of $E$ leaves 
$M \cap \Omega_{\rm alg}$ invariant.\end{proof}

We assemble these observations to finish the proof of \fullref{zeromain}. 

\begin{proof}[Proof of \fullref{zeromain}]
Since \fullref{propEregular} already proves
\fullref{zeromain} for the case of regular values, 
we only need to discuss the case when $a$ is a singular value. 
In this case, we use \fullref{econnected} and \fullref{denseon} and an argument exactly
analogous to that given for \fullref{secondmain} at the end of
\fullref{secregular}.\end{proof}

\bibliographystyle{gtart}
\bibliography{link}

\end{document}